# Output-input stability and minimum-phase nonlinear systems

Daniel Liberzon*
Coordinated Science Lab.
Univ. of Illinois at Urbana-Champaign
Urbana, IL 61801

A. Stephen Morse†
Dept. of Electrical Engineering
Yale University
New Haven, CT 06520

Eduardo D. Sontag‡
Dept. of Mathematics
Rutgers University
New Brunswick, NJ 08903

**Abstract**

This paper introduces and studies the notion of output-input stability, which represents a variant of the minimum-phase property for general smooth nonlinear control systems. The definition of output-input stability does not rely on a particular choice of coordinates in which the system takes a normal form or on the computation of zero dynamics. In the spirit of the "input-to-state stability" philosophy, it requires the state and the input of the system to be bounded by a suitable function of the output and derivatives of the output, modulo a decaying term depending on initial conditions. The class of output-input stable systems thus defined includes all affine systems in global normal form whose internal dynamics are input-to-state stable and also all left-invertible linear systems whose transmission zeros have negative real parts. As an application, we explain how the new concept enables one to develop a natural extension to nonlinear systems of a basic result from linear adaptive control.

## I Introduction

A continuous-time linear single-input/single-output (SISO) system is said to be *minimum-phase* if the numerator polynomial of its transfer function has all its zeros in the open left half of the complex plane. This property can be given a simple interpretation that involves the *relative degree* of the system, which equals the difference between the degrees of the denominator and the numerator of the transfer function. Namely, if a linear system of relative degree $r$ is minimum-phase, then the "inverse" system, driven by the $r$-th derivative of the output of the original system, is stable. For left-invertible, multi-input/multi-output (MIMO) systems, in place of the zeros of the numerator one appeals to the so-called *transmission zeros* [10].

The notion of a minimum-phase system is of great significance in many areas of linear system analysis and design. In particular, it has played an important role in parameter adaptive control. A basic example is provided by the "certainty equivalence output stabilization theorem" [11], which says that when a certainty equivalence, output stabilizing adaptive controller is applied to a minimum-phase linear system, the closed-loop system is detectable through the tuning error. In essence, this result serves as a justification for the certainty equivalence approach to adaptive control of minimum-phase linear systems.

For nonlinear systems that are affine in controls, a major contribution of Byrnes and Isidori [3] was to define the minimum-phase property in terms of the new concept of *zero dynamics*. The zero dynamics are the internal dynamics of the system under the action of an input that holds the output constantly at zero. The system is called *minimum-phase* if the zero dynamics are (globally) asymptotically stable. In the SISO case, a unique input capable of producing the zero output is guaranteed to exist if the system has a uniform relative degree, which is now defined to be the number of times one has to differentiate the output until the input appears. Extensions to MIMO systems are discussed in [4, 6].

*Corresponding author. Fax: (217) 244-1653. Email: liberzon@uiuc.edu. URL: http://black.csl.uiuc.edu/~liberzon
†Email: morse@cs.yale.edu. URL: http://entity.eng.yale.edu/controls. Supported by AFOSR, DARPA, and NSF
‡Email: sontag@control.rutgers.edu. URL: http://www.math.rutgers.edu/~sontag. Supported by US Air Force Grant F49620-98-1-0242

In view of the need to work with the zero dynamics, the above definition of a minimum-phase nonlinear system prompts one to look for a change of coordinates that transforms the system into a certain normal form. It has also been recognized that just asymptotic stability of the zero dynamics is sometimes insufficient for control design purposes, so that additional requirements need to be placed on the internal dynamics of the system. One such common requirement is that the internal dynamics be *input-to-state stable* with respect to the output and its derivatives up to order $r-1$, where $r$ is the relative degree (see, for instance, [13]). These remarks suggest that while the current notion of a minimum-phase nonlinear system is important and useful, it is also of interest to develop alternative (and possibly stronger) concepts which can be applied when asymptotic stability of zero dynamics is inadequate.

In this paper we introduce the notion of *output-input stability*, which does not rely on zero dynamics or normal forms and is not restricted to affine systems. Loosely speaking, we will call a system output-input stable if its state and input eventually become small when the output and derivatives of the output are small. Conceptually, the new notion relates to the existing concept of a minimum-phase nonlinear system in much the same way as input-to-state stability (ISS) relates to global asymptotic stability under zero inputs (0-GAS), modulo the duality between inputs and outputs. An important outcome of this parallelism is that the tools that have been developed for studying ISS and related concepts can be employed to study output-input stability, as will be discussed below.

It will follow from our definition that if a system has a uniform relative degree (in an appropriate sense) and is detectable through the output and its derivatives up to some order, uniformly over all inputs that produce a given output, then it is output-input stable. For SISO systems that are real analytic in controls, we will show that the converse is also true, thus arriving at a useful equivalent characterization of output-input stability (Theorem 1). We will prove that the class of output-input stable systems as defined here includes all left-invertible linear systems whose transmission zeros have negative real parts (Theorem 2) and all affine systems in global normal form with input-to-state stable internal dynamics.

Relying on a series of observations and auxiliary results deduced from the new definition, we will establish a natural nonlinear counterpart of the certainty equivalence output stabilization theorem from linear adaptive control (Theorem 3). This conceptually important and intuitively appealing result did not seem to be attainable within the boundaries of the existing theory of minimum-phase nonlinear systems. It serves to illustrate that output-input stability is a reasonable and useful extension of the notion of a minimum-phase linear system. In view of the remarks made earlier, it is probable that the new concept will find other applications in a variety of nonlinear control contexts.

The proposed definition is precisely stated in the next section. In Section III we give a somewhat non-standard definition of relative degree, which is especially suitable for subsequent developments. In Section IV we review the notions of detectability and input-to-state stability. In Section V we study the output-input stability property with the help of the concepts discussed in Sections III and IV. In Section VI we derive some useful results for cascade systems. In Section VII we present a nonlinear version of the certainty equivalence output stabilization theorem. Section VIII contains some remarks on output-input stability of input/output operators. The contributions of the paper are briefly summarized in Section IX. Examples are provided throughout the paper to illustrate the ideas.

## II Definition and preliminary remarks

We consider nonlinear control systems of the general form

$$\begin{aligned} \dot{x} &= f(x,u) \\ y &= h(x) \end{aligned} \tag{1}$$

where the state $x$ takes values in $\mathbb{R}^n$, the input $u$ takes values in $\mathbb{R}^m$, the output $y$ takes values in $\mathbb{R}^l$ (for some positive integers $n$, $m$, and $l$), and the functions $f$ and $h$ are smooth ($C^\infty$). Admissible input (or "control") signals are locally essentially bounded, Lebesgue measurable functions $u : [0, \infty) \to \mathbb{R}^m$. For every initial condition $x(0)$ and every input $u(\cdot)$, there is a maximally defined solution $x(\cdot)$ of the system (1), and the corresponding output $y(\cdot)$. Note that whenever the input function $u$ is $k-1$ times continuously differentiable, where $k$ is a positive integer, the derivatives $\dot{y}$, $\ddot{y}$, ..., $y^{(k)}$ are well defined (this issue will be addressed in more detail later).

Recall that a function $\alpha : [0, \infty) \to [0, \infty)$ is said to be of *class* $\mathcal{K}$ if it is continuous, strictly increasing, and $\alpha(0) = 0$. If $\alpha \in \mathcal{K}$ is unbounded, then it is said to be of *class* $\mathcal{K}_\infty$. A function $\beta : [0, \infty) \times [0, \infty) \to [0, \infty)$ is said to be of *class* $\mathcal{KL}$ if $\beta(\cdot, t)$ is of class $\mathcal{K}$ for each fixed $t \geq 0$ and $\beta(s, t)$ decreases to 0 as $t \to \infty$ for each fixed $s \geq 0$.

We will let $\|\cdot\|_{[a,b]}$ denote the essential supremum norm of a signal restricted to an interval $[a, b]$, i.e. $\|z\|_{[a,b]} := \operatorname{ess\,sup}\{|z(s)| : a \leq s \leq b\}$, where $|\cdot|$ is the standard Euclidean norm. When some vectors $v_1 \in \mathbb{R}^{n_1}, \ldots, v_k \in \mathbb{R}^{n_k}$ are given, we will often use the simplified notation $(v_1; \ldots; v_k)$ for the "stack" vector $(v_1^T, \ldots, v_k^T)^T \in \mathbb{R}^{n_1+\cdots+n_k}$. Given an $\mathbb{R}^l$-valued signal $z$ and a nonnegative integer $k$, we will denote by $\mathbf{z}^k$ the $\mathbb{R}^{l(k+1)}$-valued signal

$$\mathbf{z}^k := (z_1; \dot{z}_1; \ldots; z_1^{(k)}; \ldots; z_l; \dot{z}_l; \ldots; z_l^{(k)})$$

provided that the indicated derivatives exist.

We are now in position to introduce the main concept of this paper.

**Definition 1** We will call the system (1) *output-input stable* if there exist a positive integer $N$, a class $\mathcal{KL}$ function $\beta$, and a class $\mathcal{K}_\infty$ function $\gamma$ such that for every initial state $x(0)$ and every $N-1$ times continuously differentiable input $u$ the inequality

$$\left|\begin{pmatrix} u(t) \\ x(t) \end{pmatrix}\right| \leq \beta(|x(0)|, t) + \gamma(\|\mathbf{y}^N\|_{[0,t]}) \tag{2}$$

holds for all $t$ in the domain of the corresponding solution of (1). □

The inequality (2) can be interpreted in terms of two separate properties of the system. The first one is that if the output and its derivatives are small, then the input becomes small. Roughly speaking, this means that the system has a stable left inverse in the input-output sense. However, no explicit construction of such a left inverse is necessary. The second property is that if the output and its derivatives are small, then the state becomes small. This signifies that the system is (zero-state) detectable through the output and its derivatives, uniformly with respect to inputs[1]; we will call such systems "weakly uniformly 0-detectable" (see Section IV). Thus output-input stable systems form a subclass of weakly uniformly 0-detectable ones. Detectability is a state-space concept, whose attractive feature is that it can be characterized by Lyapunov-like dissipation inequalities.

In view of the bound on the magnitude of the state, a more complete name for the property introduced here would perhaps be "(differential) output-to-input-and-state stability", or DOSIS, but we choose to call it output-input stability for simplicity. The conditions imposed by Definition 1 capture intrinsic properties of the system, which are independent of a particular coordinate representation. They are consistent with the intuition provided by the concept of a minimum-phase linear system. In fact, for SISO linear systems the output-input stability property reduces precisely to the classical minimum-phase property, as we now show.

[1]In particular, this detectability property is preserved under state feedback; for comparison, recall that every minimum-phase linear system is detectable under all feedback laws.

**Example 1** Consider a stabilizable and detectable linear SISO system

$$\begin{aligned} \dot{x} &= Ax + bu \\ y &= c^T x \end{aligned} \tag{3}$$

Let $r$ be its relative degree. This means that we have $c^T b = c^T A b = \cdots = c^T A^{r-2} b = 0$ but $c^T A^{r-1} b := g \neq 0$. From the formula $y^{(r)}(t) = c^T A^r x(t) + g u(t)$ we immediately obtain

$$|u(t)| \leq \frac{|(A^r)^T c||x(t)| + |y^{(r)}(t)|}{|g|}. \tag{4}$$

Moreover, it is well known that there exists a linear change of coordinates $x \mapsto (\xi, \eta)$, where $\xi \in \mathbb{R}^r$, $\eta \in \mathbb{R}^{n-r}$, and $\xi_1 = y$, which transforms the system (3) into the normal form

$$\begin{aligned} \dot{\xi}_1 &= \xi_2 \\ \dot{\xi}_2 &= \xi_3 \\ &\cdots \\ \dot{\xi}_r &= d^T \xi + f^T \eta + gu \\ \dot{\eta} &= P\xi + Q\eta \end{aligned}$$

and (3) is minimum-phase (in the classical sense) if and only if $Q$ is a stable matrix. Stability of $Q$ is equivalent to the existence of positive constants $\lambda$ and $\mu$ such that for all initial states and all inputs $u$ we have $|\eta(t)| \leq e^{-\lambda t}|\eta(0)| + \mu \|\xi\|_{[0,t]}$ (it is also equivalent to detectability of the transformed system with extended output $\xi = \mathbf{y}^{r-1}$). Combining the last inequality with (4), we arrive at

$$|u(t)| \leq \frac{|(A^r)^T c| e^{-\lambda t}}{|g|} |x(0)| + \frac{|(A^r)^T c|(\mu + 1) + 1}{|g|} \|\mathbf{y}^r\|_{[0,t]}.$$

This yields (2) with $N = r$. On the other hand, if (2) holds, then we know that $y \equiv 0$ implies $x \to 0$, and so (3) must be minimum-phase. Thus we see that for stabilizable and detectable linear SISO systems, output-input stability is equivalent to the usual minimum-phase property. Incidentally, note that when $N = r$ in (2), the smoothness of $u$ becomes superfluous, because $y$ is automatically $r$ times (almost everywhere) differentiable for every admissible input $u$. □

The above remarks suggest that the concepts of relative degree and detectability are related to the output-input stability property. In the subsequent sections, we will develop some machinery which is needed to study this relationship, and explore to what extent the situation described in Example 1 carries over to (possibly MIMO) nonlinear systems.

We conclude this section with another motivating observation, important from the point of view of control design. It is formally expressed by the following proposition.

**Proposition 1** *Assume that the system* (1) *is output-input stable. Suppose that we are given a feedback law $u = k(x)$ with the following properties: there exist a system*

$$\dot{\xi} = F(\xi), \qquad \xi \in \mathbb{R}^q \tag{5}$$

*with a globally asymptotically stable equilibrium at the origin, a map $H : \mathbb{R}^q \to \mathbb{R}^l$ with $H(0) = 0$, and a class $\mathcal{K}_\infty$ function $\chi$ such that for every initial condition $x(0)$ for the system*

$$\dot{x} = f(x, k(x)) \tag{6}$$

*there exists an initial condition $\xi(0)$ for the system* (5) *with $|\xi(0)| \leq \chi(|x(0)|)$ for which we have $h(x(t)) = H(\xi(t))$ along the corresponding solutions. Then the origin is a globally asymptotically stable equilibrium of the closed-loop system* (6).

PROOF. Global asymptotic stability of the system (5) is expressed by the inequality

$$|\xi(t)| \leq \bar{\beta}(|\xi(0)|, t) \tag{7}$$

for some $\bar{\beta} \in \mathcal{KL}$. Under the action of the control law $u = k(x)$, the output of the system (1) satisfies (omitting time dependence) $y = H(\xi)$, $\dot{y} = \frac{\partial H}{\partial \xi} F(\xi)$, and so on. Using the standing assumptions on $F$ and $H$, it is straightforward to verify that for some class $\mathcal{K}_\infty$ function $\rho$ we have $|\mathbf{y}^N(t)| \leq \rho(|\xi(t)|)$, where $N$ is the integer that appears in Definition 1. Combined with (7), this gives

$$|\mathbf{y}^N(t)| \leq \rho(|\bar{\beta}(\chi(|x(0)|), t)|). \tag{8}$$

Using the inequality (2), global asymptotic stability of the closed-loop system can now be established by standard arguments (see Section VI for details). □

As an example, consider an affine system given by

$$\begin{aligned} \dot{x} &= f(x) + g(x)u \\ y &= h(x) \end{aligned} \tag{9}$$

Suppose that this system has a uniform relative degree $r$ (in the sense of [6]) and is output-input stable. It may or may not have a global normal form (construction of a global normal form requires additional properties besides relative degree [6]). However, relabeling $y$ as $\xi_1$, we have

$$\begin{aligned} \dot{\xi}_1 &= \xi_2 \\ \dot{\xi}_2 &= \xi_3 \\ &\cdots \\ \dot{\xi}_r &= b(x) + a(x)u \end{aligned}$$

where $a(x) := L_g L_f^{r-1} h(x)$ and $b(x) := L_f^r h(x)$. We can then apply a state feedback law which brings this to the form $\dot{\xi} = F(\xi)$, where the last system has a globally asymptotically stable equilibrium at the origin. (One possible choice is a linearizing feedback

$$u = -\frac{1}{a(x)}\big(b(x) + p_{r-1}\xi_r + \cdots + p_0\xi_1\big) = -\frac{1}{a(x)}\Big(b(x) + \sum_{i=0}^{r-1} p_i L_f^i h(x)\Big)$$

where the numbers $p_i$, $i = 0, \ldots, r-1$ are such that polynomial $s^r + p_{r-1}s^{r-1} + \cdots + p_1 s + p_0$ has all its roots in the left half-plane.) In view of Proposition 1, the entire system becomes globally asymptotically stable (here we can take $\chi$ to be the identity function). Note that this is true for *every* such feedback law, which for minimum-phase systems is generally not the case and special care needs to be taken in designing a stabilizing feedback (see [6, Section 9.2]).

# III Relative degree

## A. SISO systems

Let us first consider the case when the system (1) is SISO, i.e. when $m = l = 1$. To specify what will be meant by "relative degree", we need to introduce some notation. For each $k = 0, 1, \ldots$ define, recursively, the functions $H_k : \mathbb{R}^n \times \mathbb{R}^k \to \mathbb{R}$ by the formulas $H_0 := h$ and

$$H_{k+1}(x, u_0, \ldots, u_k) := \frac{\partial H_k}{\partial x} f(x, u_0) + \sum_{j=0}^{k-1} \frac{\partial H_k}{\partial u_j} u_{j+1} \tag{10}$$

where the arguments of $H_k$ are $(x, u_0, \ldots, u_{k-1})$. As an illustration, in the special case of the SISO affine system (9) we have $H_1(x, u_0) = L_f h + L_g h u_0$ and $H_2(x, u_0, u_1) = L_f^2 h + (L_g L_f h + L_f L_g h) u_0 + L_g^2 h u_0^2 + L_g h u_1$ (omitting the argument $x$ in the directional derivatives).

The significance of the functions $H_k$ lies in the fact that if the input $u(\cdot)$ is in $C^{k-1}$ for some positive integer $k$, then along each solution $x(\cdot)$ of (1) the corresponding output has a continuous $k$-th derivative satisfying

$$y^{(k)}(t) = H_k\big(x(t), u(t), \ldots, u^{(k-1)}(t)\big).$$

In particular, suppose that $H_k$ is independent of $u_0, \ldots, u_{k-1}$ for all $k$ less than some positive integer $r$. Then $H_r$ depends only on $x$ and $u_0$, as given by

$$H_r(x, u_0) = \frac{\partial H_{r-1}}{\partial x}(x)\, f(x, u_0)\,.$$

As an example, for the affine system (9) we have

$$H_r(x, u_0) = L_f^r h(x) \;+\; L_g L_f^{r-1} h(x)\, u_0\,. \tag{11}$$

In this case we see that for every initial condition and every input, $y^{(r-1)}$ exists and is an absolutely continuous function of time, and we have

$$y^{(r)}(t) = H_r(x(t), u(t))$$

for almost all $t$ in the domain of the corresponding solution. The converse is also true, namely, if $y^{(r-1)}$ exists and is absolutely continuous for all initial states and all inputs, then $H_k$ must be independent of $u_0, \ldots, u_{k-1}$ for all $k < r$. Indeed, if $H_k$ were a non-constant function of $u_0$ for some $k < r$, then it would take different values at $(x_0, u_{01})$ and $(x_0, u_{02})$ for some $x_0 \in \mathbb{R}^n$ and $u_{01}, u_{02} \in \mathbb{R}$. Choosing the initial state $x(0) = x_0$ and applying the input $u(t) \equiv u_{01}$ for $t \in [0, \varepsilon)$ and $u(t) \equiv u_{02}$ for $t \geq \varepsilon$, where $\varepsilon$ is small enough, would produce an output $y$ with a discontinuous $k$-th derivative, contradicting the existence and absolute continuity of $y^{(r-1)}$.

**Definition 2** We will say that a positive integer $r$ is the (uniform) *relative degree* of the system (1) if the following two conditions hold:

1. For each $k < r$, the function $H_k$ is independent of $u_0, \ldots, u_{k-1}$.
2. There exist two class $\mathcal{K}_\infty$ functions $\rho_1$ and $\rho_2$ such that
$$|u_0| \leq \rho_1(|x|) + \rho_2\left(|H_r(x, u_0)|\right) \tag{12}$$
for all $x \in \mathbb{R}^n$ and all $u_0 \in \mathbb{R}$. □

If there exists such an integer $r$, then it is unique. This can be deduced from the following simple observation, which will also be needed later.

**Remark 1** If properties 1 and 2 in Definition 2 hold for some positive integer $r$, then there does not exist an $x_0$ such that $H_r(x_0, u_0) = C$ for all $u_0$, where $C$ is some constant. (Indeed, otherwise we would have $|u_0| \leq \rho_1(x_0) + \rho_2(|C|)$ for all $u_0$, a contradiction.) □

We conclude that if some $r$ satisfies Definition 2, then $H_r$ cannot be independent of $u_0$ by virtue of Remark 1, hence property 1 in Definition 2 cannot hold for any $\bar{r} > r$. In view of the previous remarks, $r$ is the relative degree of (1) if and only if for some functions $\rho_1, \rho_2 \in \mathcal{K}_\infty$, for every initial condition, and every input, $y^{(r-1)}$ exists and is absolutely continuous (hence $y^{(r)}$ exists almost everywhere) and the inequality

$$|u(t)| \leq \rho_1(|x(t)|) + \rho_2(|y^{(r)}(t)|) \tag{13}$$

holds for almost all $t$ (to see why (13) implies (12), simply apply an arbitrary constant control $u(t) \equiv u_0$). We next show that for affine systems the above definition is consistent with the usual one.

**Proposition 2** *Consider the affine system* (9)*. A positive integer $r$ is the relative degree of* (9) *in the sense of Definition 2 if and only if $L_g L_f^k h(x) = 0$ for all $x$ and all integers $k < r-1$, $L_g L_f^{r-1} h(x) \neq 0$ for all $x$, and $L_f^r h(0) = 0$.*

Proof. Suppose that $r$ is the relative degree of (9) in the sense of Definition 2. Applying (11) repeatedly, we see that property 1 in Definition 2 implies $L_g L_f^k h \equiv 0$ for all $k < r-1$. Moreover, if there were some $x_0 \in \mathbb{R}^n$ with $L_g L_f^{r-1} h(x_0) = 0$, then we would have $H_r(x_0, u_0) = L_f^r h(x_0)$ which is independent of $u_0$, contrary to Remark 1. Now, setting $\bar{u}_0 := -L_f^r h(0)/L_g L_f^{r-1} h(0)$, we have $H_r(0, \bar{u}_0) = 0$, and so (12) implies that $\bar{u}_0$ must be 0 hence $L_f^r h(0) = 0$.

Conversely, suppose that $r$ satisfies the properties in the statement of the proposition. Then property 1 in Definition 2 is verified using (11). Moreover, $H_r(x, u_0)$ takes the form $b(x) + a(x)u_0$, where $a(x) := L_g L_f^{r-1} h(x)$ and $b(x) := L_f^r h(x)$. Since $a(x) \neq 0$ for all $x$ and $b(0) = 0$, we have

$$\frac{1}{|a(x)|} \leq \bar{\rho}_1(|x|) + \frac{1}{|a(0)|}, \qquad |b(x)| \leq \bar{\rho}_2(|x|)$$

for some class $\mathcal{K}_\infty$ functions $\bar{\rho}_1$ and $\bar{\rho}_2$. If follows that

$$|u_0| \leq \frac{1}{|a(x)|}\left(|H_r(x, u_0)| + |b(x)|\right) \leq \left(\bar{\rho}_1(|x|) + \frac{1}{|a(0)|}\right)\left(|H_r(x, u_0)| + \bar{\rho}_2(|x|)\right)$$

from which property 2 in Definition 2 can be easily shown to hold. □

Proposition 2 implies, in particular, that for the SISO affine system (9) with $f(0) = 0$ Definition 2 reduces to the standard definition of uniform relative degree as given, e.g., in the book by Isidori [6] (simply note that $f(0) = 0$ implies $L_f^r h(0) = 0$). Of course, the definition of relative degree proposed here is not restricted to affine systems. As a simple example, the system $\dot{y} = u^2$ has relative degree 1 according to Definition 2. This case is also covered by the definition of relative degree for not necessarily affine systems given in [12, p. 417]. In general, however, our definition is more restrictive; for example, the system $\dot{y} = \arctan u$ would have relative degree 1 in the context of [12], but the bound (12) does not hold.

**Remark 2** The relative degree cannot exceed the dimension of the system, i.e. if $r$ exists, then we must have $r \leq n$. For the affine system (9), this is an immediate consequence of Lemma 4.1.1 of [6] or Corollary 5.3.8 of [16], which imply that if $r$ is the relative degree, at least locally at some $x_0 \in \mathbb{R}^n$, then the $r$ row vectors $dh(x)$, $dL_f h(x)$, $\ldots$, $dL_f^{r-1} h(x)$ are linearly independent in a neighborhood of $x_0$ in $\mathbb{R}^n$ (so the functions $h(x)$, $L_f h(x)$, $\ldots$, $L_f^{r-1} h(x)$ qualify as a partial set of new local coordinates). For the general system (1), we can consider the associated affine system

$$\begin{aligned} \dot{x} &= f(x, z) \\ \dot{z} &= u \\ y &= h(x) \end{aligned} \tag{14}$$

of dimension $n+1$. If $r$ is the relative degree of (1), then there exists a point $(x_0, z_0) \in \mathbb{R}^{n+1}$ at which (14) has relative degree $r+1$ (in the sense of [6]). The results just mentioned then imply that the $r+1$ functions $H_0(x)$, $\ldots$, $H_r(x, z)$ can serve as independent coordinates in a neighborhood of $(x_0, z_0)$ in $\mathbb{R}^{n+1}$, hence $r \leq n$. □

To prove one of our main results (Theorem 1 in Section V) we will need the following characterization of relative degree, which also has intrinsic interest.

**Proposition 3** *A positive integer $r$ is the relative degree of* (1) *if and only if the following three conditions are satisfied:*

1. *For each $k < r$, the function $H_k$ is independent of $u_0, \dots, u_{k-1}$.*
2. *For each compact set $\mathcal{X} \subset \mathbb{R}^n$ and each positive constant $K$, there exists a number $M$ such that $|H_r(x, u_0)| \geq K$ whenever $x \in \mathcal{X}$ and $|u_0| \geq M$.*
3. *$H_r(0, u_0) \neq 0$ for all $u_0 \neq 0$.*

This in turn requires the following lemma, which is a simple exercise on $\mathcal{K}_\infty$ functions. Its proof is included for completeness.

**Lemma 1** *If a continuous function $H : \mathbb{R}^p \to [0, \infty)$, where $p$ is a positive integer, is such that $H(z) \to \infty$ when $|z| \to \infty$, and $H(z) > 0$ when $z \neq 0$, then there exists a class $\mathcal{K}_\infty$ function $\rho$ such that $|z| \leq \rho(H(z))$ for all $z$.*

Proof of Lemma 1. Consider the function $\gamma : [0, \infty) \to [0, \infty)$ given by $\gamma(s) := \min_{|z| \geq s} H(z)$. This function is well defined (because $H$ is radially unbounded so the minimum is taken over a compact set), continuous, positive definite, nondecreasing, and unbounded, and we have $\gamma(|z|) \leq H(z)$ for all $z$. Then one can find a function $\alpha \in \mathcal{K}_\infty$ such that $\alpha \leq \gamma$ (for arbitrary $s_0 \leq s_1 \leq s_2$ such that $\gamma$ is increasing on $[s_0, s_1]$ and constant on $[s_1, s_2]$, let $\alpha(s) := (\gamma(s) + \gamma(s_0))/2$ for $s \in [s_0, s_1]$, and on $[s_1, s_2]$ let $\alpha$ be linear with $\alpha(s_2) = \gamma(s_2)$). We have $\alpha(|z|) \leq H(z)$, hence $|z| \leq \rho(H(z))$ with $\rho := \alpha^{-1}$. □

Proof of Proposition 3. Condition 1 in the statement of the proposition exactly matches property 1 in Definition 2. If property 2 in Definition 2 holds, then condition 2 of the proposition is satisfied with $M(\mathcal{X}, K) = \rho_1(\max_{x \in \mathcal{X}} |x|) + \rho_2(K)$. It is clear that property 2 in Definition 2 implies condition 3 of the proposition. It remains to show that conditions 2 and 3 of the proposition imply property 2 in Definition 2. Let $z := (x, u_0)$, and consider the function $H(z) := |x| + |H_r(x, u_0)|$. We claim that this function satisfies the hypotheses of Lemma 1. Indeed, take an arbitrary $K > 0$. Condition 2 of the proposition implies that there exists an $M$ such that $|H_r(x, u_0)| > K$ whenever $|x| \leq K$ and $|u_0| > M$. It follows that $H(z) > K$ if $|z| > K + M$. Therefore, $H$ is radially unbounded. In view of condition 3 of the proposition, clearly $H(z) > 0$ if $z \neq 0$. Thus we can apply Lemma 1, which guarantees the existence of a function $\rho \in \mathcal{K}_\infty$ such that

$$\left| \begin{pmatrix} x \\ u_0 \end{pmatrix} \right| \leq \rho\big(|x| + |H_r(x, u_0)|\big) \leq \rho(2|x|) + \rho\big(2|H_r(x, u_0)|\big).$$

From this (12) follows with $\rho_1(s) = \rho_2(s) := \rho(2s)$. □

## B. MIMO systems

The above concept extends in a straightforward fashion to the case when the system (1) is MIMO, i.e. when $m$ and $l$ are not necessarily equal to 1. For each $i \in \{1, \dots, l\}$, let $H_0^i$ be the $i$-th component of $h$, and define the functions $H_k^i : \mathbb{R}^n \times (\mathbb{R}^m)^k \to \mathbb{R}$, $k = 1, 2, \dots$ recursively by the formula (10) with $H_k^i, H_{k+1}^i$ instead of $H_k, H_{k+1}$. We will say that a set of positive integers $\{r_1, \dots, r_l\}$ is a (uniform) *relative degree* of the system (1) if the following two conditions hold:

1. For each $i \in \{1, \dots, l\}$ and each $k < r_i$, the function $H_k^i$ is independent of $u_0, \dots, u_{k-1}$.

2. There exist two class $\mathcal{K}_\infty$ functions $\rho_1$ and $\rho_2$ such that

$$|u_0| \leq \rho_1(|x|) + \rho_2\Big(\Big|\big(H^1_{r_1}(x,u_0);\ldots;H^l_{r_l}(x,u_0)\big)\Big|\Big)$$

for all $x \in \mathbb{R}^n$ and all $u_0 \in \mathbb{R}^m$.

Similarly to the SISO case, $\{r_1,\ldots,r_l\}$ is the relative degree of (1) if and only if for some functions $\rho_1,\rho_2 \in \mathcal{K}_\infty$, for every $i \in \{1,\ldots,l\}$, every initial condition, and every input, $y_i^{(r_i-1)}$ exists and is absolutely continuous (hence $y_i^{(r_i)}$ exists almost everywhere) and the inequality

$$|u(t)| \leq \rho_1(|x(t)|) + \rho_2\Big(\Big|\big(y_1^{(r_1)}(t);\ldots;y_l^{(r_l)}(t)\big)\Big|\Big)$$

holds for almost all $t$. Using the same type of argument as in the proof of Proposition 2, one can show that for the square $(m = l)$ affine system

$$\begin{aligned} \dot{x} &= f(x) + \sum_{i=1}^{m} g_i(x)u_i \\ y_1 &= h_1(x) \\ &\cdots \\ y_m &= h_m(x) \end{aligned}$$

with $f(0) = 0$ this reduces to the definition of uniform vector relative degree given in [6], which says that we must have $L_{g_j}L_f^k h_i(x) = 0$ for all $x$ and all integers $1 \leq j \leq m$, $1 \leq i \leq m$, $k < r_i - 1$, and that the $m \times m$ matrix $A(x)$ defined by

$$A_{ij}(x) := L_{g_j}L_f^{r_i-1}h_i(x)$$

must be nonsingular for all $x$. More generally, for $m \leq l$ the corresponding $l \times m$ matrix must be left-invertible (i.e. of rank $m$) for all $x$. Note that for MIMO systems the relative degree is not necessarily unique (example: $\dot{y}_1 = u_1^2 + u_2^2$, $\dot{y}_2 = y_2$; $\{1,k\}$ is a relative degree for every $k > 0$). Proposition 3 carries over to the MIMO setting subject to an obvious change of notation, but this will not be needed in the sequel.

# IV Detectability and related notions

Consider a general system of the form

$$\dot{x} = f(x,u).$$

We recall from [15] that this system is called *input-to-state stable* (ISS) if there exist some functions $\beta \in \mathcal{KL}$ and $\gamma \in \mathcal{K}_\infty$ such that for every initial state $x(0)$ and every input $u$ the corresponding solution satisfies the inequality

$$|x(t)| \leq \beta(|x(0)|,t) + \gamma(\|u\|_{[0,t]})$$

for all $t \geq 0$. Intuitively, this means that the state eventually becomes small when the input is small.

Given a system with both inputs and outputs

$$\begin{aligned} \dot{x} &= f(x,u) \\ y &= h(x,u) \end{aligned} \tag{15}$$

we will say that it is *0-detectable* if there exist some functions $\beta \in \mathcal{KL}$ and $\gamma_1,\gamma_2 \in \mathcal{K}_\infty$ such for every $x(0)$ and every $u$ the corresponding solution satisfies the inequality

$$|x(t)| \leq \beta(|x(0)|,t) + \gamma_1(\|u\|_{[0,t]}) + \gamma_2(\|y\|_{[0,t]})$$

as long as it exists. In particular, a system without inputs given by

$$\dot{x} = f(x)$$
$$y = h(x)$$

will be called *0-detectable* if there exist some functions $\beta \in \mathcal{KL}$ and $\gamma \in \mathcal{K}_\infty$ such that for every initial state $x(0)$ the corresponding solution satisfies the following inequality as long as it exists:

$$|x(t)| \le \beta(|x(0)|, t) + \gamma(\|y\|_{[0,t]}). \tag{16}$$

These concepts were studied in [17] under the names of *input/output-to-state stability* and *output-to-state stability*, respectively. In this paper we use the term "0-detectability" as a shorthand.

Let us call the system (15) *uniformly 0-detectable* if there exist some functions $\beta \in \mathcal{KL}$ and $\gamma \in \mathcal{K}_\infty$ such that for every initial state $x(0)$ and every input $u$ the inequality (16) holds along the corresponding solution. As the name suggests, uniform 0-detectability amounts to 0-detectability that is uniform with respect to inputs. This property was called *uniform output-to-state stability* in [7] (see also [8]) and *strong detectability* in [5].

When working with output derivatives, it is helpful to introduce an auxiliary system whose output contains derivatives of the output of the original system (restricting admissible inputs if necessary to ensure that these derivatives are well defined). We first describe this construction for SISO systems. Consider the system (1) with $m = l = 1$. Take a nonnegative integer $k$. Restricting the input $u$ to be in $C^{k-1}$, we can consider the *k-output extension* of (1):

$$\begin{aligned} \dot{x} &= f(x, u) \\ \mathbf{y}^k &= \mathbf{h}_k(x, u, \dots, u^{(k-1)}) \end{aligned} \tag{17}$$

where

$$\mathbf{h}_k(x, u, \dots, u^{(k-1)}) := \big(H_0(x); H_1(x, u); \dots; H_k(x, u, \dots, u^{(k-1)})\big)$$

is the new output map (here we are using the notation of Section III). That is, we redefine the output of the system to be $\mathbf{y}^k$. Of course, for $k = 0$ we recover the original system. Note that the differentiability assumption on $u$ can be relaxed if the function $H_k$ is known to be independent of $u_0, \dots, u_{k-1}$ for all $k$ between 0 and some positive integer (cf. Section III). We will view (17) as a system whose input consists of $u$ and all derivatives of $u$ that appear as arguments in the output map $\mathbf{h}_k$. With some abuse of terminology, we will apply to such systems the definitions of 0-detectability and uniform 0-detectability given earlier. Let us call the system (1) *weakly 0-detectable of order k* if its $k$-output extension (17) is 0-detectable. Also, let us call the system (1) *weakly uniformly 0-detectable of order k* if its $k$-output extension (17) is uniformly 0-detectable.

Now consider the general MIMO case, i.e. take the system (1) with $u \in \mathbb{R}^m$ and $y \in \mathbb{R}^l$. For arbitrary nonnegative numbers $k_1, \dots, k_l$, we can redefine the output of the system to be

$$\mathbf{y}^{k_1,\dots,k_l} := (y_1; \dot{y}_1; \dots; y_1^{(k_1)}; \dots; y_l; \dot{y}_l; \dots; y_l^{(k_l)})$$

restricting the input $u$ to be sufficiently smooth so that the indicated derivatives exist. This amounts to considering the system

$$\begin{aligned} \dot{x} &= f(x, u) \\ \mathbf{y}^{k_1,\dots,k_l} &= \mathbf{h}_{k_1,\dots,k_l}(x, u, \dots, u^{(j-1)}) \end{aligned} \tag{18}$$

where $j := \max_{1\le i\le l} k_i$ and

$$\mathbf{h}_{k_1,\dots,k_l}(x,u,\dots,u^{(j-1)}) := \begin{pmatrix} \Big(H_0^1(x), H_1^1(x,u), \dots, H_{k_1}^1(x,u,\dots,u^{(k_1-1)})\Big)^T \\ \cdots \\ \Big(H_0^l(x), H_1^l(x,u), \dots, H_{k_l}^l(x,u,\dots,u^{(k_l-1)})\Big)^T \end{pmatrix}$$

(in the notation of Section III). Similarly to the SISO case, we will call this system the $\{k_1,\dots,k_l\}$-*output extension* of (1). We will say that the system (1) is *weakly 0-detectable of order* $\{k_1,\dots,k_l\}$ if its $\{k_1,\dots,k_l\}$-output extension (18) is 0-detectable. We will also say that the system (1) is *weakly uniformly 0-detectable of order* $\{k_1,\dots,k_l\}$ if its $\{k_1,\dots,k_l\}$-output extension (18) is uniformly 0-detectable.

Another definition, which will be needed in Section VII, is the following one (introduced in [18]). The system (15) is said to be *input-to-output stable* if there exist some functions $\beta \in \mathcal{KL}$ and $\gamma \in \mathcal{K}_\infty$ such that for every $x(0)$ and every $u$ the following inequality holds along the corresponding solution:

$$|y(t)| \le \beta(|x(0)|, t) + \gamma(\|u\|_{[0,t]}).$$

Finally, we remark that since all the systems under consideration are time-invariant, the same properties would result if we used an arbitrary initial time $t_0$ instead of 0 in the above inequalities (changing the second argument of $\beta$ from $t$ to $t - t_0$ accordingly). This fact will be used implicitly in the proofs in Section VI.

# V Output-input stability

## A. SISO systems

We first study the SISO case, represented by the system (1) with $m = l = 1$. Our main result in this section is the following characterization of output-input stability for SISO systems.

**Theorem 1** *1. Suppose that the system* (1) *has a relative degree $r$ in the sense of Definition 2 and is weakly uniformly 0-detectable of order $k$, for some $k$. Then* (1) *is output-input stable in the sense of Definition 1, with $N = \max\{r, k\}$.*

*2. Suppose that the system* (1) *is output-input stable in the sense of Definition 1. Then it is weakly uniformly 0-detectable of order $N$.*

*3. Suppose that the system* (1) *is output-input stable in the sense of Definition 1, that the function $f(x,\cdot)$ is real analytic in $u$ for each fixed $x$, and that $f(0,0) = 0$ and $h(0) = 0$. Then* (1) *has a relative degree $r \le N$ in the sense of Definition 2.*

The theorem implies that for systems with relative degree, output-input stability is equivalent to weak uniform 0-detectability of order $k$ for some $k \le N$, and that for systems satisfying the additional assumptions on $f$ and $h$ stated in part 3, output-input stability is equivalent to the existence of a relative degree $r \le N$ plus weak uniform 0-detectability of order $k$ for some $k \le N$.

Proof. Part 1. Since $r$ is the relative degree, the inequality (13) holds with $\rho_1, \rho_2 \in \mathcal{K}_\infty$. Suppose that for some $k$ the system (1) is weakly uniformly 0-detectable of order $k$. This can be expressed as

$$|x(t)| \le \bar\beta(|x(0)|, t) + \bar\gamma(\|\mathbf{y}^k\|_{[0,t]})$$

where $\bar\beta \in \mathcal{KL}$ and $\bar\gamma \in \mathcal{K}_\infty$ (see Section IV). Combining this with (13) and using the simple fact that for every class $\mathcal{K}$ function $\rho$ and arbitrary numbers $s_1, s_2 \ge 0$ one has $\rho(s_1 + s_2) \le \rho(2s_1) + \rho(2s_2)$, we arrive at the inequality (2) with $N := \max\{r, k\}$, $\beta(s,t) := \rho_1(2\bar{\beta}(s,t)) + \bar{\beta}(s,t)$, and $\gamma(s) := \rho_1(2\bar{\gamma}(s)) + \rho_2(s) + \bar{\gamma}(s)$. Thus (1) is output-input stable as needed.

Part 2. Follows immediately from the definitions.

Part 3. Since the system (1) is output-input stable, we know in particular that for some positive integer $N$ and some functions $\beta \in \mathcal{KL}$ and $\gamma \in \mathcal{K}_\infty$ the inequality

$$|u(t)| \le \beta(|x(0)|, t) + \gamma(\|\mathbf{y}^N\|_{[0,t]}) \tag{19}$$

holds along solutions of (1) for all smooth inputs. The function $H_k$ cannot be independent of $u_0, \ldots, u_{k-1}$ for all $k = 0, \ldots, N$. Otherwise, letting $x(0) = 0$ and applying an arbitrary constant input $u_0$, we would deduce from (19) that $|u_0| \le c := \beta(0,t) + \gamma\big(\big|(H_0(0); \ldots; H_N(0))\big|\big)$, a contradiction. Thus the integer

$$r := \max\{k < N : H_k \text{ is independent of } u_0, \ldots, u_{k-1}\} + 1 \tag{20}$$

is well defined. Condition 1 of Proposition 3 holds with this $r$. For every input $u$ we have:

$$\begin{aligned}
\dot{y} &= H_1(x) \\
&\cdots \\
y^{(r-1)} &= H_{r-1}(x) \\
y^{(r)} &= \frac{\partial H_{r-1}}{\partial x}(x) f(x,u) = H_r(x,u) \\
y^{(r+1)} &= \frac{\partial H_r}{\partial x}(x,u) f(x,u) + \frac{\partial H_r}{\partial u_0}(x,u)\dot{u} =: G_1(x,u) + \frac{\partial H_r}{\partial u_0}(x,u)\dot{u} \\
y^{(r+2)} &=: G_2(x,u,\dot{u}) + \frac{\partial H_r}{\partial u_0}(x,u)\ddot{u} \\
&\cdots \\
y^{(N)} &=: G_{N-r}(x,u,\ldots,u^{(N-r-1)}) + \frac{\partial H_r}{\partial u_0}(x,u) u^{(N-r)}
\end{aligned} \tag{21}$$

The following fact will be useful.

**Lemma 2** *If* (19) *holds and $r$ is defined by* (20)*, then there cannot exist a bounded sequence $\{x_j\}$ in $\mathbb{R}^n$, a sequence $\{w_j\}$ in $\mathbb{R}$ with $\lim_{j\to\infty} |w_j| = \infty$, and a positive constant $K$ such that for all $j$ we have $|H_r(x_j, w_j)| < K$ and $\frac{\partial H_r}{\partial u_0}(x_j, w_j) \neq 0$.*

Proof. Suppose that there exist sequences $\{x_j\}$ and $\{w_j\}$ and a positive constant $K$ with the properties indicated in the statement of the lemma. Fix an arbitrary positive integer $j$. Consider the initial state $x(0) = x_j$, and pick a smooth (e.g., polynomial) input function $u_j(\cdot)$ with $u_j(0) = w_j$ whose derivatives at $t = 0$ are specified recursively by the equations

$$u_j^{(i)}(0) = -\frac{G_i(x_j, u_j(0), \ldots, u_j^{(i-1)}(0))}{\frac{\partial H_r}{\partial u_0}(x_j, u_j(0))}, \qquad i = 1, \ldots, N - r. \tag{22}$$

In view of (21), we will then have

$$\begin{aligned}
&\dot{y}(0) = H_1(x_j), \ \ldots, \ y^{(r-1)}(0) = H_{r-1}(x_j), \\
&|y^{(r)}(0)| = |H_r(x_j, w_j)| < K, \quad y^{(r+1)}(0) = \cdots = y^{(N)}(0) = 0.
\end{aligned}$$

Therefore, if such an input $u_j$ is applied and if $\varepsilon$ is an arbitrary fixed positive number, then there exists a sufficiently small time $T_j$ such that for all $t \in [0, T_j]$ the following inequalities hold:

$$|x(t)| < |x_j| + \varepsilon, \quad |\dot{y}(t)| < |H_1(x_j)| + \varepsilon, \ \ldots, \ |y^{(r-1)}(t)| < |H_{r-1}(x_j)| + \varepsilon,$$
$$|y^{(r)}(t)| < K + \varepsilon, \quad |y^{(r+1)}(t)| < \varepsilon, \ \ldots, \ |y^{(N)}(t)| < \varepsilon.$$

Repeating this construction for all $j$, we obtain a sequence of trajectories of (1) along which $|x(t)|$, $|\dot{y}(t)|$, $\ldots$, $|y^{(N)}(t)|$ are uniformly bounded for small $t$, whereas $|u_j(t)|$ is unbounded for small $t$ and large $j$. We arrive at a contradiction with (19), and the proof of the lemma is complete. □

Let us denote by $\Theta$ the set of all $x \in \mathbb{R}^n$ such that $H_r(x, \cdot)$ is a constant function. The set $\Theta$ is closed (because if for a sequence $\{x_j\}$ converging to some $x$ we have $\frac{\partial H_r}{\partial u_0}(x_j, u_0) = 0$ for all $u_0$ and all $j$, then $\frac{\partial H_r}{\partial u_0}(x, u_0) = 0$ for all $u_0$).

**Lemma 3** *Suppose that $f(x, \cdot)$ is real analytic in $u$ for each fixed $x$, that (19) holds, and that $r$ is defined by (20). If $\bar{x} \in \mathbb{R}^n$ is such that for some $K > 0$ we have $|H_r(\bar{x}, u)| < K$ for all $u$, then $\bar{x}$ is in the interior of $\Theta$.*

Proof. Take an arbitrary sequence $\{v_j\}$ in $\mathbb{R}$ with $\lim_{j\to\infty} |v_j| = \infty$. By hypothesis, $|H_r(\bar{x}, v_j)| < K$ for all $j$. By continuity, for each $j$ there exist a neighborhood $B_j(\bar{x})$ of $\bar{x}$ in $\mathbb{R}^n$ and a positive number $\delta_j$ such that $|H_r(x, u_0)| < K$ for all $x \in B_j(\bar{x})$ and all $u_0 \in (v_j - \delta_j, v_j + \delta_j)$. Moreover, the neighborhoods $B_j(\bar{x})$, $j = 1, 2, \ldots$ can be chosen to be nested, i.e. $B_i(\bar{x}) \subset B_j(\bar{x})$ whenever $i > j$, and the sequence $\{\delta_j\}$ can be chosen to be nonincreasing. Now, suppose that $\bar{x}$ is not in the interior of $\Theta$. Fix an arbitrary $j$. We have $B_j(\bar{x}) \not\subset \Theta$. Take an arbitrary $x_j \in B_j(\bar{x}) \setminus \Theta$. Then $\frac{\partial H_r}{\partial u_0}(x_j, \cdot)$ cannot be identically zero on the interval $(v_j - \delta_j, v_j + \delta_j)$, by virtue of real analyticity of $H_r(x_j, \cdot)$ which follows from that of $f(x_j, \cdot)$. Thus we can find a $w_j \in (v_j - \delta_j, v_j + \delta_j)$ such that $\frac{\partial H_r}{\partial u_0}(x_j, w_j) \neq 0$. This construction can be carried out for all $j$. Since $B_i(\bar{x}) \subset B_j(\bar{x})$ when $i > j$, the points $x_j$, $j = 1, 2, \ldots$ can be chosen in such a way that $|x_j|$ is uniformly bounded for all $j$. Moreover, we have $|w_j| \geq |v_j| - |\delta_j| \geq |v_j| - |\delta_1| \to \infty$ as $j \to \infty$. In view of Lemma 2, we arrive at a contradiction with (19), which proves the lemma. □

**Corollary 1** *Suppose that $f(x, \cdot)$ is real analytic in $u$ for each fixed $x$. If (19) holds and $r$ is defined by (20), then the set $\Theta$ is open.*

Proof. By definition of $\Theta$, every $\bar{x} \in \Theta$ satisfies the condition in the statement of Lemma 3, hence it lies in the interior of $\Theta$. □

**Corollary 2** *Suppose that $f(x, \cdot)$ is real analytic in $u$ for each fixed $x$. If (19) holds and $r$ is defined by (20), then the set $\Theta$ is empty.*

Proof. We know that $\Theta$ is closed. We also know that $\Theta$ is open (Corollary 1). Moreover, $\Theta \neq \mathbb{R}^n$ by virtue of (20). Being a closed and open proper subset of $\mathbb{R}^n$, the set $\Theta$ must be empty. □

Our goal is to show that conditions 2 and 3 of Proposition 3 hold, which would imply that $r$ is the relative degree of (1). We break this into two separate statements.

**Lemma 4** *Suppose that $f(x, \cdot)$ is real analytic in $u$ for each fixed $x$. If (19) holds and $r$ is defined by (20), then for each compact set $\mathcal{X} \subset \mathbb{R}^n$ and each positive constant $K$ there exists a number $M$ such that $|H_r(x, u_0)| \geq K$ whenever $x \in \mathcal{X}$ and $|u_0| \geq M$.*

Proof. Suppose that, contrary to the statement of the lemma, there exist a compact subset $\mathcal{X}$ of $\mathbb{R}^n$, a positive constant $K$, a sequence $\{x_j\}$ in $\mathcal{X}$, and a sequence $\{v_j\}$ in $\mathbb{R}$ with $\lim_{j\to\infty} |v_j| = \infty$ such that $|H_r(x_j, v_j)| < K$ for all $j$. By continuity, we can find a nonincreasing sequence of positive numbers $\{\delta_j\}$

such that $|H_r(x_j, u_0)| < K$ for all $j$ and all $u_0 \in (v_j - \delta_j, v_j + \delta_j)$. Fix an arbitrary $j$. Since $\Theta$ is empty by Corollary 2 and $H_r(x_j, \cdot)$ is real analytic, $\frac{\partial H_r}{\partial u_0}(x_j, \cdot)$ cannot vanish identically on the interval $(v_j - \delta_j, v_j + \delta_j)$. Thus we can find a $w_j \in (v_j - \delta_j, v_j + \delta_j)$ such that $\frac{\partial H_r}{\partial u_0}(x_j, w_j) \neq 0$. Repeat this construction for all $j$. Lemma 2 applies again, yielding a contradiction with (19), and the proof of the lemma is complete. □

**Lemma 5** *Suppose that $f(x, \cdot)$ is real analytic in $u$ for each fixed $x$ and that we have $f(0, 0) = 0$ and $h(0) = 0$. If (19) holds and $r$ is defined by (20), then $H_r(0, u_0) \neq 0$ for all $u_0 \neq 0$.*

PROOF. Suppose that $H_r(0, \bar{u}_0) = 0$ for some $\bar{u}_0 \neq 0$. We know from Corollary 2 that the set $\Theta$ is empty. Thus by real analyticity $\frac{\partial H_r}{\partial u_0}(0, \cdot)$ cannot vanish identically on any open neighborhood of $\bar{u}_0$. This implies that there exists a sequence $\{v_j\}$ converging to $\bar{u}_0$ such that $\frac{\partial H_r}{\partial u_0}(0, v_j) \neq 0 \ \forall j$ (if $\frac{\partial H_r}{\partial u_0}(0, \bar{u}_0) \neq 0$, simply let $v_j \equiv \bar{u}_0$). Choose an arbitrary $j$. Take the initial state to be $x(0) = 0$. Pick a smooth (e.g., polynomial) input function $u_j(\cdot)$ such that $u_j(0) = v_j$ and the equations (22) hold with 0 in place of $x_j$. From (21) we immediately see that $y^{(r+1)}(0) = \cdots = y^{(N)}(0) = 0$. Since $h(0) = 0$, we have $y(0) = 0$. We also know that $f(0, 0) = 0$, which implies that $H_1(0) = \cdots = H_{r-1}(0) = 0$. It follows that $\dot{y}(0) = \cdots = y^{(r-1)}(0) = 0$. We conclude that if the input $u_j$ is applied, then for every $\varepsilon > 0$ there exists a sufficiently small time $T_j$ such that for all $t \in [0, T_j]$ the following inequalities hold:

$$
\begin{aligned}
&|x(t)| < \varepsilon, \quad |\dot{y}(t)| < \varepsilon, \ \ldots, \ |y^{(r-1)}(t)| < \varepsilon, \\
&|y^{(r)}(t)| < |H_r(0, v_j)| + \varepsilon, \quad |y^{(r+1)}(t)| < \varepsilon, \ \ldots, \ |y^{(N)}(t)| < \varepsilon.
\end{aligned}
$$

Carrying out the above construction for all $j$ and noting that $\lim_{j \to \infty} H_r(0, v_j) = H_r(0, \bar{u}_0) = 0$, we see that $y(t), \dot{y}(t), \ldots, y^{(N)}(t)$ become arbitrarily small for small $t$ as $j \to \infty$. On the other hand, $u_j(0) \to \bar{u}_0 \neq 0$, so $u_j(t)$ is bounded away from 0 for small $t$ and large $j$. This is a contradiction with (19), which proves the lemma. □

We have shown that the integer $r$ defined by (20) satisfies all three conditions of Proposition 3, thus $r$ is the relative degree of the system (1). This proves part 3, and the proof of Theorem 1 is complete. □

As an illustration, consider the affine system (9) with $f(0) = 0$. Its right-hand side is obviously real analytic in $u$. Reconstructing the above proof for this case, we find that $r$ is the smallest integer for which $L_g L_f^{r-1} h(x)$ is not identically zero on $\mathbb{R}^n$, and $\Theta = \{x : L_g L_f^{r-1} h(x) = 0\}$. If the system is output-input stable, then Corollary 2 implies that $\Theta$ must be empty, which means that $r$ is the relative degree (see Proposition 2). Since the hypothesis $h(0) = 0$ is only used in Lemma 5, it is not needed in this case.

We will be especially interested in systems that are covered by part 1 of Theorem 1 with $k = r - 1$. We give such systems a separate name to emphasize the relationship with the existing terminology (which will be explained in the next example). Note that $\mathbf{y}^{r-1}$ is a function of the state $x$ only: $\mathbf{y}^{r-1} = \mathbf{h}_{r-1}(x)$; no differentiability assumptions need to be placed on $u$.

**Definition 3** Let us call the system (1) *strongly minimum-phase* if it has a relative degree $r$ and is weakly uniformly 0-detectable of order $r - 1$. □

**Example 2** Consider an affine system in global normal form

$$
\begin{aligned}
\dot{\xi}_1 &= \xi_2 \\
\dot{\xi}_2 &= \xi_3 \\
&\cdots \\
\dot{\xi}_r &= b(\xi, \eta) + a(\xi, \eta) u \\
\dot{\eta} &= q(\xi, \eta) \\
y &= \xi_1
\end{aligned}
\tag{23}
$$

where $\xi := (\xi_1; \dots; \xi_r)$, $b(0,0) = 0$, and $a(\xi, \eta) \neq 0 \ \forall \xi, \eta$ (so that $r$ is the uniform relative degree by Proposition 2). This system is usually called minimum-phase if the *zero dynamics* $\dot{\eta} = q(0, \eta)$ have an asymptotically stable equilibrium at $\eta = 0$ (see [3]). Since $\mathbf{y}^{r-1} = \xi$, the above definition of the strong minimum-phase property in this case demands that the equation for $\eta$ in (23), which represents the internal dynamics, be input-to-state stable (ISS) with respect to $\xi$ (more precisely, with respect to all possible signals $\xi$ that can be generated by the $\xi$-subsystem). This is in general a stronger condition than just asymptotic stability of the zero dynamics; however, in the linear case the two properties are equivalent (and both amount to saying that all zeros of the transfer function must have negative real parts). As we already mentioned, the ISS assumption has been imposed on the internal dynamics of the system in various contexts associated with control design (see, e.g., [13]). □

**Remark 3** The bound (13), which is a consequence of our definition of relative degree, does not necessarily imply that one can explicitly solve for $u$ in terms of $x$ and $y^{(r)}$ ("flatness"); example: $\dot{y} = u^2$. However, this is possible for some systems, in particular, it can always be done for the affine system (9). In this case, we can express $u$ as a function of $x$ and $\bar{y} := p(D)y$, where $p$ is an arbitrary stable polynomial of degree $r$ and $D := \frac{d}{dt}$. Substituting this expression for $u$, we obtain an "inverse" system, driven by $\bar{y}$. If the system (1) is strongly minimum-phase, then it is not hard to show that this inverse system will be ISS with respect to $\bar{y}$. For example, consider the system in global normal form (23). Take a stable polynomial $p(s) = s^r + p_{r-1}s^{r-1} + \dots + p_1 s + p_0$, and rewrite the equation for $\dot{\xi}_r$ as $\dot{\xi}_r = -p_{r-1}\xi_r - \dots - p_1\xi_2 - p_0\xi_1 + \bar{y}$. Then the $r$-dimensional subsystem that describes the evolution of $\xi$ is easily seen to be a stable linear system driven by $\bar{y}$, hence it is ISS with respect to $\bar{y}$. If the system $\dot{\eta} = q(\xi, \eta)$ is ISS with respect to $\xi$ (recall that this is a consequence of the strong minimum-phase property) then the overall system is indeed ISS with respect to $\bar{y}$, because a cascade of two ISS systems is ISS. □

The results of [8, 17] imply that the system (1) is weakly uniformly 0-detectable of order $r-1$ if there exists a smooth, positive definite, radially unbounded function $V : \mathbb{R}^n \to \mathbb{R}$ that satisfies

$$\nabla V(x) f(x,u) \leq -\alpha(|x|) + \chi(|\mathbf{y}^{r-1}|) \qquad \forall x, u \tag{24}$$

for some functions $\alpha, \chi \in \mathcal{K}_\infty$. This Lyapunov-like dissipation inequality can be used to check the strong minimum-phase property, once the relative degree of the system is known. We summarize this observation in the following statement.

**Proposition 4** *Suppose that the system* (1) *has a relative degree $r$ and that for some smooth, positive definite, radially unbounded function $V : \mathbb{R}^n \to \mathbb{R}$ and class $\mathcal{K}_\infty$ functions $\alpha, \chi$ the inequality* (24) *holds. Then* (1) *is strongly minimum-phase.*

In fact, the inequality (24) provides a necessary and sufficient condition for weak uniform 0-detectability of order $r-1$ if controls take values in a compact set [8]. Unfortunately, this condition is only sufficient and not necessary if the control set is unbounded. For example, consider the integrator $\dot{x} = u$, $y = x$. It is obviously uniformly 0-detectable (here $r = 1$), but for every smooth positive definite $V : \mathbb{R} \to \mathbb{R}$ and every $x$ in the nonempty set $\{x \in \mathbb{R} : V'(x) \neq 0\}$ the quantity $V'(x)u$ can be made arbitrarily large by a suitable choice of $u$.

Of course, the inequality (24) with $N$ instead of $r-1$ can be applied to check weak uniform 0-detectability of order $N$, as long as $y^{(N)}$ is well defined. A similar recipe for finding the relative degree with the help of Lyapunov-like functions does not seem to exist. However, for most systems of interest it is not difficult to verify the existence of relative degree directly by using Definition 2 or Proposition 3.

**Example 3** We now give an example of an output-input stable system that is not strongly minimum-phase. Consider the system

$$\begin{aligned} \dot{x}_1 &= u \\ \dot{x}_2 &= -x_2 + u^2 \\ y &= x_1 \end{aligned} \tag{25}$$

It has relative degree 1. From the equation $\dot{x}_2 = -x_2 + \dot{y}^2$, which is ISS with respect to $\dot{y}$, we see that the system (25) is weakly uniformly 0-detectable of order 1. Therefore, (25) is output-input stable (with $N = 1$) by virtue of part 1 of Theorem 1. Now let us show that this system is not strongly minimum-phase, i.e. it is not uniformly 0-detectable (with respect to the original output $y$). It is enough to find a solution trajectory along which $x_1$ converges to zero while $x_2$ does not converge to zero. Take the initial state to be 0, and apply the following input: $u(t) \equiv 1$ for $0 \le t < 1$, $u(t) \equiv -1$ for $1 \le t < 2$, $u(t) \equiv 1$ for $2 \le t < 2\frac{1}{2}$, $u(t) \equiv -1$ for $2\frac{1}{2} \le t < 3$, $u(t) \equiv 1$ for $3 \le t < 3\frac{1}{4}$, and so on. Then $x_1 \to 0$, whereas for $x_2$ we have $\dot{x}_2 = -x_2 + 1$ so $x_2 \to 1$. □

The situation in the above example is dual to the one described in [1], where it is shown that ISS with respect to inputs and their derivatives is in general not equivalent to the usual ISS.

## B. MIMO systems

We now turn to the general case of the system (1) with $u \in \mathbb{R}^m$ and $y \in \mathbb{R}^l$. The next result is readily obtained by the same arguments as those employed in the proof of Theorem 1. We will see below that the nontrivial part of Theorem 1, which states that under suitable assumptions output-input stability implies the existence of a relative degree, does not hold for MIMO systems. We will also explain why this is an advantage, rather than a drawback, of our definition.

**Proposition 5** *Suppose that the system* (1) *has a relative degree* $\{r_1, \dots, r_l\}$. *Then it is output-input stable, with* $N \ge \max\{r_1, \dots, r_l\}$, *if and only if it is weakly uniformly 0-detectable of order* $\{k_1, \dots, k_l\}$ *for some* $k_1, \dots, k_l \le N$.

Of special interest are systems satisfying the condition of Proposition 5 with $k_i = r_i - 1$, $i = 1, \dots, l$. Accordingly, let us call the system (1) *strongly minimum-phase* if it has a relative degree $\{r_1, \dots, r_l\}$ and is weakly uniformly 0-detectable of order $\{r_1 - 1, \dots, r_l - 1\}$. Note that $\mathbf{y}^{r_1-1,\dots,r_l-1}$ is a function of $x$ only, as given by $\mathbf{y}^{r_1-1,\dots,r_l-1} = \mathbf{h}_{r_1-1,\dots,r_l-1}(x)$, and no differentiability assumptions need to be placed on $u$. We thus obtain a generalization to MIMO systems of the strong minimum-phase property introduced in the previous subsection. It has a similar interpretation in terms of input-to-state stability of the internal dynamics for systems in global normal form, and admits an analogous Lyapunov-like sufficient condition.

However, we remark that for MIMO systems the existence of a relative degree is quite a restrictive assumption. For example, linear systems with relative degree form a rather special subclass of those linear systems for which the minimum-phase property (in its classical sense) is well defined.[2] Fortunately, Definition 1 does not have the shortcoming of applying only to systems with relative degree, as illustrated by the following example.

[2]Basically, the reason for this is that the relative degree of a MIMO system can be lost or gained as a result of a linear coordinate transformation in the output space, while the minimum-phase property is invariant under such transformations.

**Example 4** The system

$$\begin{aligned}
\dot{x}_1 &= u_1 \\
\dot{x}_2 &= x_3 + u_1^2 \\
\dot{x}_3 &= u_2 \\
\dot{x}_4 &= -x_4 + x_1^3 \\
y &= (x_1; x_2)
\end{aligned}$$

does not have a relative degree. This system is output-input stable, as can be seen from the formulas $|u_2| = |\ddot{y}_2 - 2\dot{y}_1\ddot{y}_1| \leq |\ddot{y}_2| + \dot{y}_1^2 + \ddot{y}_1^2$, $|x_3| \leq |\dot{y}_2| + \dot{y}_1^2$, and the fact that the equation for $x_4$ is ISS with respect to $x_1$. □

The above example is to be contrasted with the next one.

**Example 5** The system

$$\begin{aligned}
\dot{x}_1 &= u_1 \\
\dot{x}_2 &= x_3 + x_2 u_2 \\
\dot{x}_3 &= u_2 \\
\dot{x}_4 &= -x_4 + x_1^3 \\
y &= (x_1; x_2)
\end{aligned}$$

does not have a relative degree. The zero dynamics of this system are $\dot{x}_4 = -x_4$, and the input that produces them is identically zero. However, it can be shown by the same kind of argument as the one used in the proof of Theorem 1 that arbitrarily large $u_2$ can lead to arbitrarily small $x$, $y$, and derivatives of $y$. Therefore, this system is not output-input stable. □

We will now establish an important feature of Definition 1, namely, that for linear MIMO systems it reduces exactly to the classical definition of the minimum-phase property. We will make use of some concepts and results from the linear geometric control theory (see [10, 19]). Consider the linear system

$$\begin{aligned}
\dot{x} &= Ax + Bu \\
y &= Cx
\end{aligned} \tag{26}$$

with $x \in \mathbb{R}^n$, $u \in \mathbb{R}^m$, and $y \in \mathbb{R}^l$. We assume that (26) is left-invertible. Recall that a subspace $\mathcal{V}$ of $\mathbb{R}^n$ is called $(A, B)$-*invariant* if there exists an $m \times n$ matrix $F$ such that $(A + BF)\mathcal{V} \subset \mathcal{V}$. Denote by $\mathcal{I}$ the family of all $(A, B)$-invariant subspaces that are contained in $\ker C$. Then $\mathcal{I}$ has a unique largest member (with respect to inclusion) which we denote by $\mathcal{S}$. The eigenvalues of the restriction of $A + BF$ to $\mathcal{S}$ are the same for all $F$ such that $(A + BF)\mathcal{S} \subset \mathcal{S}$. These eigenvalues are called the *transmission zeros* of the system (26). Left-invertible linear systems are usually called minimum-phase if all their transmission zeros have negative real parts. Our goal is to prove that this is equivalent to the output-input stability property in the sense of Definition 1.

As a direct consequence of Silverman's "structure algorithm" [14], there exists an $m$-vector $\mathbf{y}$, whose components are linear combinations of the components of $y$ and their derivatives, which satisfies

$$\dot{\mathbf{y}} = Lx + Mu$$

where the matrix $M$ is nonsingular. From this we immediately obtain

$$|u| \leq |M^{-1}L||x| + |M^{-1}||\dot{\mathbf{y}}|.$$

An $(n-m)$-vector $z$ of complementary coordinates can be chosen whose dynamics are independent of $u$, as given by $\dot{z} = Nx$. Consider the feedback matrix $F := -M^{-1}L$. Then $\mathcal{S}$ is precisely the largest $(A+BF)$-invariant subspace in $\ker C$, and thus equals the unobservable subspace of $(C, A+BF)$, i.e. $\mathcal{S} = \bigcap_{i=1}^{n} \ker C(A+BF)^{i-1}$. It follows that there exists a linear change of coordinates $x \mapsto (\xi, \eta)$ such that $\mathcal{S} = \{(\xi,\eta) : \xi = 0\}$, and the components of $\xi$ are linear combinations of the components of $y$ and their derivatives. In these coordinates the system (26) takes the form

$$\begin{aligned} \dot{\xi} &= R\xi + GT\eta + Gu \\ \dot{\eta} &= P\xi + Q\eta \end{aligned}$$

Its transmission zeros are the eigenvalues of $Q$. Reasoning exactly as in Example 1, we see that $Q$ is a stable matrix if and only if the system is output-input stable. We summarize this as follows.

**Theorem 2** *A left-invertible linear system is output-input stable if and only if all its transmission zeros have negative real parts.*

# VI Cascade results

The purpose of this section is to investigate how the output-input stability property behaves under series connections of several subsystems. We will prove two lemmas. The first one says that the cascade of a 0-detectable system with an output-input stable system is weakly 0-detectable (i.e. 0-detectable through the output and its derivatives), which is a result of independent interest. The second lemma is a direct generalization of the first one, and will be needed to prove the adaptive control result of Section VII. To simplify the presentation and to obtain the sharpest results possible, we restrict our attention to SISO strongly minimum-phase systems (see Definition 3). The same proof techniques apply readily to SISO output-input stable systems, and in particular to systems that have a relative degree and satisfy the hypotheses of part 1 of Theorem 1 with $k \geq r$. However, for $k > r$ the conclusions become weaker. Generalizations to MIMO systems are also straightforward, subject to similar limitations.

Suppose that we are given two systems:

$$\Sigma_1 : \quad \begin{aligned} \dot{x}_1 &= f_1(x_1, u_1) \\ y_1 &= h_1(x_1) \end{aligned} \tag{27}$$

and

$$\Sigma_2 : \quad \begin{aligned} \dot{x}_2 &= f_2(x_2, u_2) \\ y_2 &= h_2(x_2) \end{aligned} \tag{28}$$

Upon setting $u_2 = y_1$, we obtain a cascade system with input $u_1$ and output $y_2$, which we denote by $\Sigma_c$ (see Figure 1).

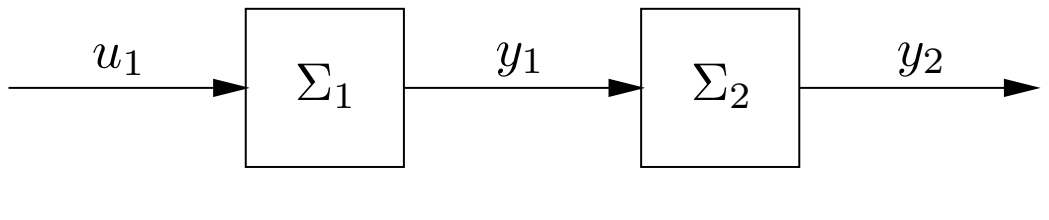

Figure 1: The cascade system

Assume that $\Sigma_2$ has a relative degree $r$. Consider the $r$-output extension of $\Sigma_c$. As explained in Section III, the extended output map has the form

$$\mathbf{h}(x_1, x_2) = \big(h_2(x_2), H_1(x_2), \ldots, H_{r-1}(x_2), H_r(x_2, h_1(x_1))\big)$$

which is independent of $u_1$. In particular, no differentiability assumptions on the inputs are needed. Thus the $r$-output extension of $\Sigma_c$ is a system with input $u_1$ and output $\mathbf{y}_2^r$. The following result says that this system is 0-detectable.

**Lemma 6** *If $\Sigma_1$ is 0-detectable and $\Sigma_2$ is strongly minimum-phase, then the cascade system $\Sigma_c$ is weakly 0-detectable of order $r$.*

PROOF. In the proofs of this lemma and the next one, $\beta$ with various subscripts will be used to denote class $\mathcal{KL}$ functions, and $\gamma$ and $\rho$ with various subscripts will be used to denote class $\mathcal{K}_\infty$ functions. For $t \geq t_0 \geq 0$, the 0-detectability of $\Sigma_1$ can be expressed by the inequality

$$|x_1(t)| \leq \beta_1(|x_1(t_0)|, t - t_0) + \gamma_0(\|u_1\|_{[t_0,t]}) + \gamma_1(\|y_1\|_{[t_0,t]})$$

while the strong minimum-phase property of $\Sigma_2$ leads to the inequalities

$$|x_2(t)| \leq \beta_2(|x_2(t_0)|, t - t_0) + \gamma_2(\|\mathbf{y}_2^{r-1}\|_{[t_0,t]}) \tag{29}$$

and

$$|y_1(t)| \leq \rho_1(|x_2(t)|) + \rho_2(|y_2^{(r)}(t)|).$$

Since the system $\Sigma_c$ has a cascade structure, we employ the trick of breaking a time interval under consideration into several parts in order to derive the result (as done in [15]). Straightforward but lengthy calculations yield

$$\begin{aligned} |x_1(t)| &\leq \beta_1(|x_1(t/2)|, t/2) + \gamma_0(\|u_1\|_{[t/2,t]}) + \gamma_1(\|y_1\|_{[t/2,t]}) \\ &\leq \bar\beta_1(|x_1(0)|, t) + \bar\beta_2(|x_2(0)|, t) + \bar\gamma_0(\|u_1\|_{[0,t]}) + \bar\gamma_1(\|\mathbf{y}_2^r\|_{[0,t]}) \end{aligned}$$

where

$$\begin{aligned} \bar\beta_1(s,t) &:= \beta_1(3\beta_1(s,t/2), t/2) \\ \bar\beta_2(s,t) &:= \beta_1(9\gamma_1(3\rho_1(2\beta_2(s,0))), t/2) + \gamma_1(3\rho_1(2\beta_2(s,t/2))) \\ \bar\gamma_0(s) &:= \gamma_0(s) + \beta_1(3\gamma_0(s), 0) \\ \bar\gamma_1(s) &:= \beta_1(9\gamma_1(3\rho_1(2\gamma_2(s))), 0) + \beta_1(9\gamma_1(3\rho_2(s)), 0) + \gamma_1(3\rho_1(2\gamma_2(s))) + \gamma_1(3\rho_2(s)) \end{aligned}$$

Combining this with (29), we arrive at the desired result. □

Lemma 6 states that the cascade system is 0-detectable through $y_2$ and derivatives of $y_2$ up to order $r$, where $r$ is the relative degree of $\Sigma_2$. It is in general not true that the cascade system is 0-detectable through the original output $y_2$ only. We support this claim by constructing an example of a 0-detectable system which, when followed by an integrator, fails to remain 0-detectable.

**Example 6** Assume given a smooth function $\varphi : [0, \infty) \to [0, \infty)$ with the following properties:

1. $\varphi$ is in $L^1$.
2. $\int_0^a \varphi^2(s)\, ds \geq a^3/4$ for all $a \geq 0$ (in particular, $\varphi$ is not in $L^2$).
3. $\varphi(0) = 0$.
4. $\varphi'(0^+) = 1$, $\varphi^{(i)}(0^+) = 0$ for $i > 1$.

(It is an elementary exercise to obtain such a $\varphi$. Start with a piecewise constant function that takes the value $k^5$ on an interval of length $1/k^7$ around each positive integer $k$, and is zero elsewhere. Since $\sum_{k=1}^{\infty} k^5 k^{-7} < \infty$, this function is in $L^1$; we also have $\sum_{k=1}^{j} k^{10} k^{-7} = 1 + \ldots + j^3$. Then approximate this function by a smooth one and modify it in a neighborhood of zero to achieve desired behavior there. Properties 2 and 4 are not conflicting, because property 4 means that $\varphi(s) \approx s$ near 0, and $\int_0^a s^2 ds = a^3/3 > a^3/4$. Property 4 is only needed to ensure the smoothness of the system $\Sigma_1$ to be constructed next.)

Consider the system $\Sigma_1$ given by (27), with $x_1 \in \mathbb{R}$, $f_1(x_1) \equiv 1$ (no input $u_1$), $h_1(x_1) := -x_1$ for $x_1 < 0$ and $h_1(x_1) := \varphi(x_1)$ for $x_1 \geq 0$. To verify that this system is 0-detectable, take the function $V$ defined for $x_1 < 0$ by $V(x_1) := x_1^2$ and for $x_1 \geq 0$ by

$$V(x_1) := \int_0^{x_1} (2\varphi^2(s) - s^2)\, ds\,.$$

The function $V$ is radially unbounded and positive definite (because for $x_1 > 0$ we have $V(x_1) \geq x_1^3/6$ using property 2 of $\varphi$). Its derivative $V'(x_1)$ equals $2x_1$ for $x_1 < 0$ and $-x_1^2 + 2\varphi^2(x_1)$ for $x_1 \geq 0$. Thus we have $V'(x_1)f_1(x_1) \leq -x_1^2 + 2y_1^2$ for all $x_1$. In view of the results of [17], this implies that $\Sigma_1$ is 0-detectable.

For $\Sigma_2$, we take an integrator (which is strongly minimum-phase), i.e. (28) with $x_2 \in \mathbb{R}$, $f_2(x_2, u_2) = u_2$, and $h_2(x_2) = x_2$. Then the cascade system $\Sigma_c$ has the form

$$\begin{aligned} \dot{x}_1 &= 1 \\ \dot{x}_2 &= h_1(x_1) \\ y_2 &= x_2 \end{aligned}$$

With initial state 0 we have $y_2(t) = \int_0^t \varphi(s) ds$, which is bounded in light of property 1 of $\varphi$, while $x_1(t) \to \infty$. $\square$

Next, suppose that the system $\Sigma_1$ has another output $y_3 = h_3(x_1)$. Letting $u_2 = y_1$ as before, and defining the output $y_4 := y_3 - y_2$, we obtain a cascade-feedforward system $\Sigma_{cf}$ with input $u_1$ and output $y_4$, shown in Figure 2.

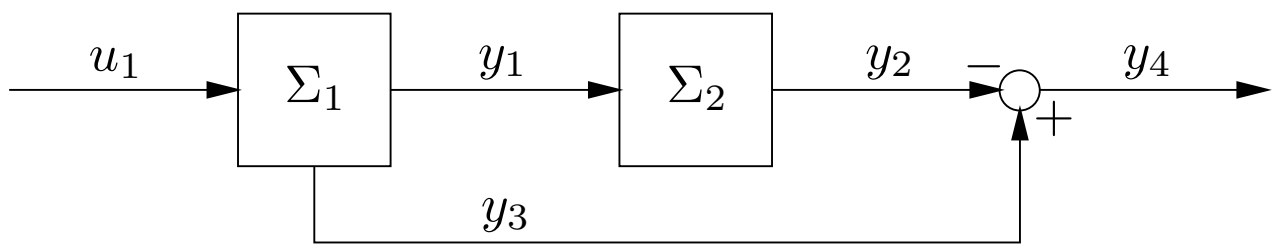

Figure 2: The cascade-feedforward system

Assume that the input $u_1$ is in $C^{r-1}$, where $r$ is the relative degree of $\Sigma_2$ as before. We can then consider the system $\widetilde{\Sigma}_1$ whose input is $\mathbf{u}_1^{r-1}$ and whose output is $\mathbf{y}_3^r$. Indeed, as explained in Section III, for each $i \in \{1, \ldots, r\}$ the $i$-th derivative of $y_3$ exists and can be written as $y^{(i)}(t) = H_i(x(t), u_1(t), \ldots, u_1^{i-1}(t))$ for a suitable function $H_i$. Moreover, since $y_2$ is $r$ times differentiable almost everywhere, we can consider the $r$-output extension of $\Sigma_{cf}$, whose input is $\mathbf{u}_1^{r-1}$ and whose output is $\mathbf{y}_4^r$. The next result says that this last system is 0-detectable.

**Lemma 7** *Suppose that $\Sigma_1$ is 0-detectable (with respect to its input $u_1$ and both its outputs, $y_1$ and $y_3$), $\Sigma_2$ is strongly minimum-phase, and the system $\widetilde{\Sigma}_1$ with input $\mathbf{u}_1^{r-1}$ and output $\mathbf{y}_3^r$ is input-to-output stable. Then the cascade-feedforward system $\Sigma_{cf}$ is weakly 0-detectable of order $r$.*

PROOF. For $t \geq t_0 \geq 0$, the hypotheses of the lemma lead to the following inequalities:

$$
\begin{aligned}
|x_1(t)| &\leq \beta_1(|x_1(t_0)|, t-t_0) + \gamma_0(\|u_1\|_{[t_0,t]}) + \gamma_1(\|y_1\|_{[t_0,t]}) + \gamma_3(\|y_3\|_{[t_0,t]}) \\
|x_2(t)| &\leq \beta_2(|x_2(t_0)|, t-t_0) + \gamma_2(\|\mathbf{y}_2^{r-1}\|_{[t_0,t]}) \\
|y_1(t)| &\leq \rho_1(|x_2(t)|) + \rho_2(|y_2^{(r)}(t)|) \\
|\mathbf{y}_3^r(t)| &\leq \beta_3(|x_1(t_0)|, t-t_0) + \gamma_4(\|\mathbf{u}_1^{r-1}\|_{[t_0,t]})
\end{aligned}
$$

We have

$$
|x_2(t)| \leq \hat{\beta}_1(|x_1(0)|, t) + \hat{\beta}_2(|x_2(0)|, t) + \hat{\gamma}_0(\|\mathbf{u}_1^{r-1}\|_{[0,t]}) + \hat{\gamma}_4(\|\mathbf{y}_4^r\|_{[0,t]})
$$

where

$$
\begin{aligned}
\hat{\beta}_1(s,t) &:= \beta_2(6\gamma_2(4\beta_3(s,0)), t/2) + \gamma_2(4\beta_3(s,t/2)) \\
\hat{\beta}_2(s,t) &:= \beta_2(3\beta_2(s,t/2), t/2) \\
\hat{\gamma}_0(s) &:= \beta_2(6\gamma_2(4\gamma_4(s)), 0) + \gamma_2(4\gamma_4(s)) \\
\hat{\gamma}_4(s) &:= \gamma_2(2s) + \beta_2(3\gamma_2(2s), 0)
\end{aligned}
$$

Similarly,

$$
|x_1(t)| \leq \tilde{\beta}_1(|x_1(0)|, t) + \tilde{\beta}_2(|x_2(0)|, t) + \tilde{\gamma}_0(\|\mathbf{u}_1^{r-1}\|_{[0,t]}) + \tilde{\gamma}_4(\|\mathbf{y}_4^r\|_{[0,t]})
$$

where, for example,

$$
\tilde{\beta}_2(s,t) := \beta_1(12\gamma_1(3\rho_1(2\beta_2(s,0))), t/2) + \gamma_1(4\rho_1(4\beta_2(2\beta_2(s,t/4), t/4))).
$$

Combining the two estimates, we obtain the desired result. □

# VII Adaptive control

In this section we describe a framework for adaptive control of uncertain nonlinear systems, in which the concept of an output-input stable system introduced in this paper turns out to be useful (enabling one to achieve what did not seem possible with the existing definition of a minimum-phase nonlinear system).

## A. Set-up and motivation

Let $\mathbb{P}$ be an unknown process, with dynamics of the form

$$
\begin{aligned}
\dot{x}_{\mathbb{P}} &= f_{\mathbb{P}}(x_{\mathbb{P}}, u) \\
y &= h_{\mathbb{P}}(x_{\mathbb{P}})
\end{aligned}
$$

where $x_{\mathbb{P}} \in \mathbb{R}^n$ is the state, $u \in \mathbb{R}$ is the control input, and $y \in \mathbb{R}$ is the measured output (we assume that $\mathbb{P}$ is SISO just to simplify the notation; the generalization to the MIMO setting is straightforward). Assume that $\mathbb{P}$ is a member of some family of systems $\bigcup_{p\in\mathcal{P}} \mathcal{F}_p$, where $\mathcal{P}$ is an index set. For each $p \in \mathcal{P}$, the subfamily $\mathcal{F}_p$ can be viewed as consisting of a *nominal process model* $\mathbb{M}_p$ together with a collection of its "perturbed" versions. The present discussion is quite general and does not depend on any special structure of $\mathcal{F}_p$.

Consider the following family of controllers, parameterized by $p$ taking values in $\mathcal{P}$:

$$
\begin{aligned}
\dot{x}_{\mathbb{C}} &= f_{\mathbb{C}}(x_{\mathbb{C}}, y, p) \\
u_p &= h_{\mathbb{C}}(x_{\mathbb{C}}, p)
\end{aligned}
$$

For every fixed $p \in \mathcal{P}$, we denote the corresponding controller by $\mathbb{C}_p$. One can think of $\mathbb{C}_p$ as a *candidate controller*, which would be used to control the process $\mathbb{P}$ if this process were known to be a member of $\mathcal{F}_p$.

We assume that on-line controller selection is carried out with the help of some estimation procedure. This is facilitated by a dynamical system $\mathbb{E}$ called the *multi-estimator*, which takes the form

$$\begin{aligned} \dot{x}_{\mathbb{E}} &= f_{\mathbb{E}}(x_{\mathbb{E}}, y, u) \\ y_p &= h_p(x_{\mathbb{E}}), \qquad p \in \mathcal{P} \end{aligned}$$

The signals $y_p, p \in \mathcal{P}$ are used to define the *estimation errors*

$$e_p := y_p - y, \qquad p \in \mathcal{P}.$$

One usually designs the multi-estimator in such a way that $e_p$ converges to zero in the case when the unknown process coincides with the $p$-th nominal process model $\mathbb{M}_p$ and there are no disturbances or noise.

Most of the standard adaptive algorithms are based on varying the index of the candidate controller in the feedback loop according to a tuning/switching law $\sigma : [0, \infty) \to \mathcal{P}$, in such a way that the corresponding estimation error $e_\sigma$ is maintained small in some sense. The underlying principle behind such a strategy is known as *certainty equivalence*. Intuitively, the motivation here is that the nominal process model with the smallest estimation error "best" approximates the actual process, and thus the candidate controller associated with that model can be expected to do the best job of controlling the process. To justify this paradigm, one must be able to ensure that the smallness of the estimation error implies the smallness of the state of the closed-loop system. Thus we see that a crucial desired property of this system is 0-detectability through the estimation error.

To make this discussion more precise, take an arbitrary fixed $q \in \mathcal{P}$. The closed-loop system, which results when the $q$-th candidate controller $\mathbb{C}_q$ is placed in the feedback loop with the process $\mathbb{P}$ and the multi-estimator $\mathbb{E}$, is described by the equations

$$\begin{aligned} \dot{x}_{\mathbb{P}} &= f_{\mathbb{P}}(x_{\mathbb{P}}, h_{\mathbb{C}}(x_{\mathbb{C}}, q)) \\ \dot{x}_{\mathbb{E}} &= f_{\mathbb{E}}(x_{\mathbb{E}}, h_{\mathbb{P}}(x_{\mathbb{P}}), h_{\mathbb{C}}(x_{\mathbb{C}}, q)) \\ \dot{x}_{\mathbb{C}} &= f_{\mathbb{C}}(x_{\mathbb{C}}, h_{\mathbb{P}}(x_{\mathbb{P}}), q) \end{aligned} \tag{30}$$

We will take the output of this system to be the estimation error $e_q = h_q(x_{\mathbb{E}}) - h_{\mathbb{P}}(x_{\mathbb{P}})$. A glance at Figure 3 might be helpful at this point. The above remarks suggest that it is desirable to design the system (30) so as to make it 0-detectable with respect to $e_q$.

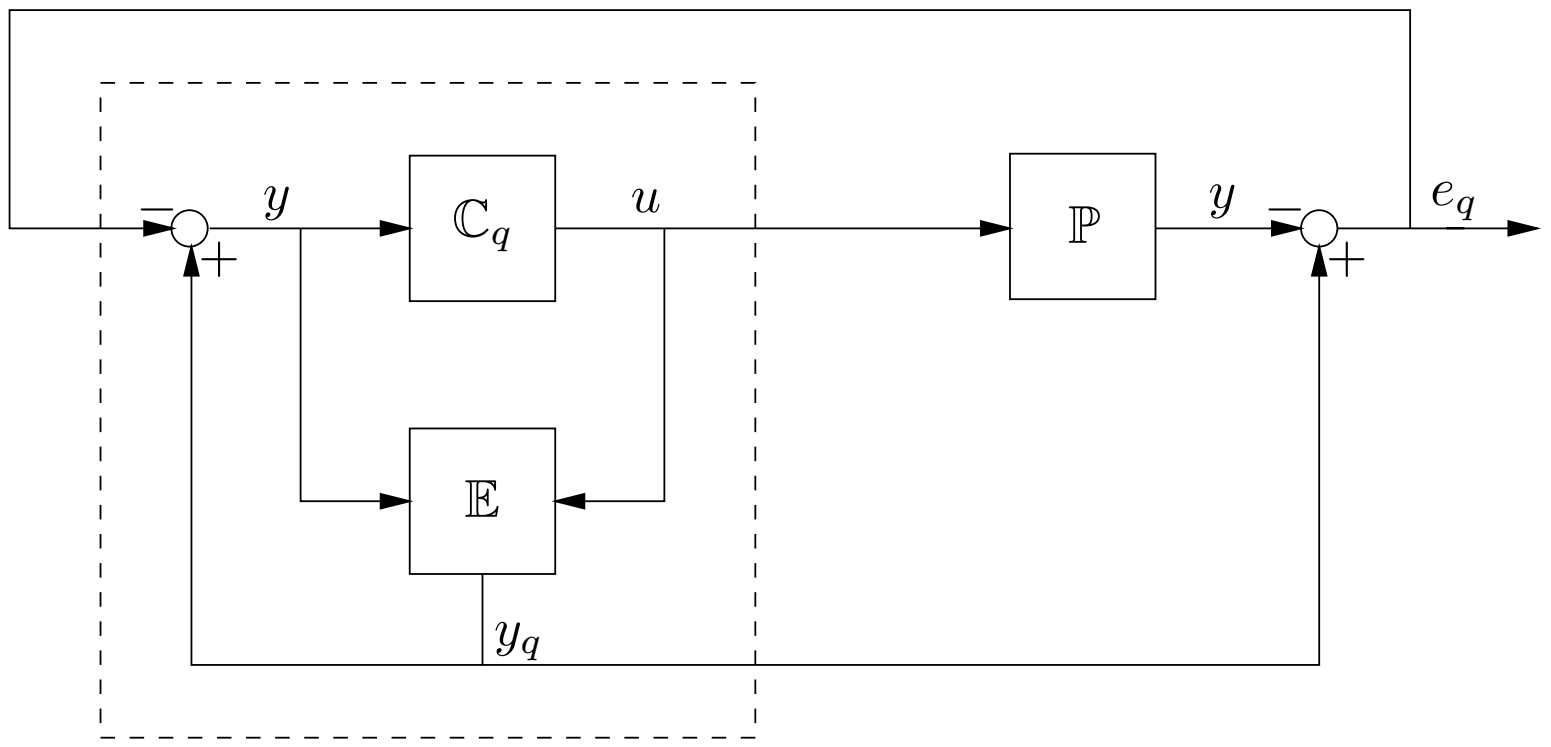

Figure 3: The closed-loop system (30)

Consider the following system, which we call the *injected system* and denote by $\mathbb{EC}_q$:

$$\begin{aligned}\dot{x}_{\mathbb{E}} &= f_{\mathbb{E}}(x_{\mathbb{E}}, h_q(x_{\mathbb{E}}) - e_q, h_{\mathbb{C}}(x_{\mathbb{C}}, q)) \\ \dot{x}_{\mathbb{C}} &= f_{\mathbb{C}}(x_{\mathbb{C}}, h_q(x_{\mathbb{E}}) - e_q, q)\end{aligned}$$

We view it as a system with input $e_q$, state $(x_{\mathbb{E}}, x_{\mathbb{C}})$, and outputs $u$ and $y_q = h_q(x_{\mathbb{E}})$. It realizes the interconnection of the $q$-th candidate controller $\mathbb{C}_q$ with the multi-estimator $\mathbb{E}$. This is the system enclosed in the dashed box in Figure 3. Basically, the choice of the candidate controllers is governed by the resulting properties of this system. This makes sense because $\mathbb{E}$ is implemented by the control designer; the validity of such an approach will become clear in view of the results that are discussed next.

It was shown in [5] that if the injected system $\mathbb{EC}_q$ is input-to-state stable (ISS) with respect to $e_q$ and the process $\mathbb{P}$ is 0-detectable, then the closed-loop system (30) is 0-detectable with respect to $e_q$. This provided a natural nonlinear extension of the Certainty Equivalence Stabilization Theorem proved for linear systems in [11]. Another relevant result from [11] is the so-called Certainty Equivalence *Output* Stabilization Theorem, which we mentioned in the Introduction. It suggests that the desired 0-detectability of the system (30) through $e_q$ should be preserved if one weakens the assumptions on the injected system $\mathbb{EC}_q$ by only requiring input-to-output stability from $e_q$ to $y_q$ instead of input-to-state stability[3], but demands that the process be output-input stable rather than 0-detectable. In what follows, we demonstrate that a result along these lines indeed holds for nonlinear systems.

## B. Main result and discussion

Assume that $\mathbb{P}$ has a known relative degree $r$. Let us redefine the input and the output of the system $\mathbb{EC}_q$ to be $\mathbf{e}_q^{r-1}$ and $\mathbf{y}_q^r$, respectively. We denote the resulting system by $\widetilde{\mathbb{EC}}_q$; its output map is obtained as explained in the previous sections. We now make the following assumptions:

1. The process $\mathbb{P}$ is strongly minimum-phase.
2. The system $\widetilde{\mathbb{EC}}_q$ is input-to-output stable.
3. The controller $\mathbb{C}_q$ is 0-detectable.
4. The multi-estimator $\mathbb{E}$ is 0-detectable.

The result given below states that the $r$-output extension of the closed-loop system (30) is 0-detectable with respect to its output $\mathbf{e}_q^r$ (the derivatives of $e_q$ exist since we only consider smooth systems; also note that there are no inputs). This result is a direct consequence of Lemma 7: one needs to apply that lemma with $\Sigma_1 = \mathbb{EC}_q$ (which is easily seen to be a 0-detectable system) and $\Sigma_2 = \mathbb{P}$.

**Theorem 3** *Under assumptions 1–4, the closed-loop system* (30) *is weakly 0-detectable of order $r$.*

The same techniques would apply readily if the process $\mathbb{P}$ is output-input stable but not necessarily strongly minimum-phase. In particular, if $\mathbb{P}$ satisfies the hypotheses of part 1 of Theorem 1 with $k \geq r$, we would conclude weak 0-detectability of order $k$ for the closed-loop system. However, for $k > r$ this is a weaker statement than that provided by Theorem 3.

Theorem 3 gives weak 0-detectability, i.e. 0-detectability through the output and derivatives of the output up to order $r$. For linear systems the distinction between 0-detectability and weak 0-detectability disappears; this is most easily seen from the well-known observability decomposition. It is possible to employ similar ideas to single out a class of nonlinear systems for which weak 0-detectability implies 0-detectability. Namely, suppose that in some coordinates the system under consideration takes the form

$$\begin{aligned}\dot{x}_1 &= f_1(x_1)\\ \dot{x}_2 &= f_2(x_1, x_2)\\ y &= h(x_1)\end{aligned} \tag{31}$$

with $f_1(0) = 0$ and $h(0) = 0$, where the subsystem $\dot{x}_1 = f_1(x_1)$, $y = h(x_1)$ is 0-detectable. This can be compared with the observability decompositions for nonlinear systems (see, e.g., [2] or [12, Chapter 3]); we replace observability with 0-detectability (neither property is weaker than the other) and require the decomposition to be global. Now, suppose that (31) is weakly 0-detectable of order $k$, where $k$ is some positive integer. Since $|\mathbf{y}^k| \le \rho(|x_1|)$, where $\rho$ is a suitable class $\mathcal{K}_\infty$ function (cf. proof of Proposition 1), we must have $|x_2(t)| \le \beta(|x_2(0)|, t) + \gamma(\|x_1\|_{[0,t]})$ for some $\beta \in \mathcal{KL}$ and $\gamma \in \mathcal{K}_\infty$. Combining this with 0-detectability of the $x_1, y$ subsystem and using the same arguments as the ones employed in Section VI for dealing with cascade systems, we see that (31) is indeed 0-detectable (through the original output $y$).

[3] A careful examination of the classical model reference adaptive control algorithm for linear systems reveals that the control law only output stabilizes the estimated model.

Therefore, if the closed-loop system (30) with output $e_q$ admits a decomposition of the kind described above, then under the assumptions of Theorem 3 this system is 0-detectable. In this case we arrive at a strengthened version of Theorem 3 which is more suitable for adaptive control applications. Of course, the existence of such a decomposition may be difficult to verify in practice. We also point out that the need to maintain smallness of the estimation error together with its first derivative is not uncommon in adaptive control (see, e.g., Chapters 5 and 6 of [9]). Precise implications of Theorem 3 for adaptive control of nonlinear systems remain to be investigated.

# VIII Output-input stability of input/output operators

It is possible to define the output-input stability property for input/output operators, without relying on state-space representations. It then turns out that if a state-space system is minimal in a suitable sense and its input/output mapping is output-input stable, then the system is output-input stable in the sense of Definition 1. Here we make some brief preliminary remarks on this; details will be pursued elsewhere.

We will consider causal operators of the form $F : \mathcal{U} \to \mathcal{Y}$, where $\mathcal{U}$ and $\mathcal{Y}$ are some spaces of time functions taking values in $\mathbb{R}^m$ and $\mathbb{R}^l$, respectively. To avoid technicalities, we can restrict these functions to be smooth almost everywhere, with domain $[0, \infty)$. What follows closely parallels the developments of [15] (see Proposition 7.1 in that paper). Let us call an operator $F$ *output-input stable* if there exist a positive integer $N$, a class $\mathcal{KL}$ function $\beta$, and a class $\mathcal{K}_\infty$ function $\gamma$ such that we have

$$|u(t)| \le \beta(\|\mathbf{y}^N\|_{[0,t_0]}, t - t_0) + \gamma(\|\mathbf{y}^N\|_{[t_0,t]}) \qquad \forall\, t \ge t_0 \ge 0 \tag{32}$$

whenever $y = F(u)$.

We think of such an operator as arising from a state-space system of the form (1), with $x(0) = 0$. Suppose that this system satisfies the following two assumptions:

1 ("observability"). There exist functions $\alpha_1, \alpha_2 \in \mathcal{K}_\infty$ such that

$$|x(t)| \le \alpha_1(\|u\|_{[t,\infty)}) + \alpha_2(\|\mathbf{y}^N\|_{[t,\infty)}) \qquad \forall\, t \ge 0.$$

This is a weak version of the strong observability property considered in [15], since the right-hand side contains derivatives of the output.

2 ("reachability"). There exists a function $\alpha_3 \in \mathcal{K}_\infty$ such that for each given $x$ it is possible to find a time $t_0$ and a control input $u$ which drives the system from state 0 at time $t = 0$ to state $x$ at time $t = t_0$ in such a way that

$$|\mathbf{y}^N(t)| \le \alpha_3(|x(t_0)|) \qquad \forall\, 0 \le t \le t_0.$$

Under appropriate conditions, this property can be derived from the strong reachability property considered in [15].

As before, the smoothness assumption on the input can be relaxed if the system has a known relative degree. It is now a simple matter to prove the following.

**Proposition 6** *If the system* (1) *satisfies assumptions 1 and 2 and the associated input/output operator is output-input stable as defined here, then* (1) *is output-input stable in the sense of Definition 1.*

# IX Conclusions

We introduced a new concept of output-input stability, which can be viewed as an ISS-like variant of the minimum-phase property for general smooth nonlinear control systems and reduces to the classical minimum-phase property for linear systems. We provided characterizations of output-input stability in terms of suitably defined notions of detectability and relative degree, the latter of which was proposed here and is of independent interest. Implications of the output-input stability property for feedback stabilization, analysis of cascade systems, and adaptive control were investigated. Other applications are likely to be found in various areas of nonlinear control theory in which the concept of a minimum-phase nonlinear system has been useful.

**Acknowledgment.** Helpful discussions with João Hespanha, Achim Ilchmann, Brian Ingalls, Alberto Isidori, and Misha Krichman are gratefully acknowledged.